\documentclass[a4paper,12pt]{article}
\usepackage[cp1251]{inputenc}        % Кодировка входного документа;
\usepackage[english,russian]{babel} % Включение русификации, русских и
\usepackage{amsmath,amsfonts,amsthm,amssymb,amsbsy,amstext,amscd,amsxtra,multicol}
\usepackage{array,mdwtab}
\usepackage[matrix,arrow,curve]{xy}
\usepackage[dvips]{graphics}
\newcommand\G{\operatorname{\mathbb G}}
\newcommand\Z{\operatorname{\mathbb Z}}
\newcommand\F{\operatorname{\mathbb F}}
\newcommand\R{\operatorname{\mathbb R}}
\renewcommand\C{\operatorname{\mathbb C}}
\newcommand\Q{\operatorname{\mathbb Q}}
\newcommand\A{\operatorname{\mathbb A}}
\newcommand\Dim{\operatorname{Dim}}

\newcommand\Det{\operatorname{Det}}
\newcommand\Frac{\operatorname{Frac}}
\newcommand\Func{\operatorname{Func}}
\newcommand\Spec{\operatorname{Spec}}
\newcommand\GL{\operatorname{GL}}
\newcommand\Pic{\operatorname{Pic}}
\newcommand\SL{\operatorname{SL}}
\newcommand\Hom{\operatorname{Hom}}
\newcommand\K{\operatorname{\mathbb K}}
\newcommand\Gal{\operatorname{Gal}}
\newcommand\Nm{\operatorname{Nm}}

\begin{document}
\thispagestyle{empty}
{\centering
\vskip2cm
{\Large {\bf Questions and remarks to the Langlands program}}\footnote{I am grateful to R. P. Langlands for very useful conversations during his visit to the Steklov Mathematical institute of the Russian Academy of Sciences (Moscow,  October 2011), to Michael Harris and Ulrich Stuhler, who  answered my sometimes too naive questions, and to \.Ilhan \.Ikeda who has read a first version of the text and has made several remarks. The work was supported by the RFBR grants  1107.0408, 1012.0486 and NSh-5139.2012.1.}
\vskip5mm
{\large A.~N.~Parshin}
\vskip2mm
({\em Uspekhi Matem. Nauk}, {\bf 67}(2012), n 3, 115-146; \\{\em Russian Mathematical Surveys}, {\bf 67}(2012), n 3, 509-539)
\vskip5mm}
%%%%%%%%%%%%%%%%%%%%%%%%%%%%%%%%%%%%%%%%%%%%%%%%%%%%%%%%%%%%
% the following few lines were added by Konstantin Besov
%%%%%%%%%%%%%%%%%%%%%%%%%%%%%%%%%%%%%%%%%%%%%%%%%%%%%%%%%%%%
\makeatletter
\def\section{\@startsection{section}{1}{-15pt}%
                           {3.5ex \@plus -1ex \@minus -.2ex}%
                           {2.3ex \@plus.2ex}%
                           {\normalfont\large\bfseries\centering\noindent}}
\def\thesection{}
\def\l@section#1#2{#1\dotfill #2}
\def\numberline#1{\par}
%\tableofcontents
{{\rightskip5mm\leftskip\rightskip\parindent0mm
%\textbf{Содержание}%
\par\smallskip\small\@starttoc{toc}\par}}
\makeatother
%%%%%%%%%%%%%%%%%%%%%%%%%%%%%%%%%%%%%%%%%%%%%%%%%%%%%%%%%%%%
\section{Introduction}
 The goal of the Langlands program is a correspondence between representations of the Galois groups  (and their generalizations or versions) and representations of reductive algebraic groups. The starting point for the construction is a field. Six types of the fields are considered: three types of local fields and three types of global fields  \cite{L3, F2}. The former ones are the following:
\vskip2mm
1) finite extensions of the field  ${\Q}_p$ of $p$-adic numbers, the field  ${\R}$ of real numbers and the field   ${\C}$ of complex numbers,

2) the fields  ${\F}_q ((t))$ of  Laurent power series, where ${\F}_q$ is the finite field of $q$ elements,

3)  the field of  Laurent power series ${\C} ((t))$.
\vskip2mm
The global fields are:
\vskip2mm
4) fields of algebraic numbers (= finite extensions of the field $\Q$ of rational numbers),

5) fields of algebraic functions in one variable with a finite constant field  (= finite extensions of the field $\F_q(x)$),

6) fields of algebraic functions in one variable with the field   $\C$ as the constant field  (= finite extensions of the field  $\C(x)$).
\vskip2mm
The local fields in the first list are completions of the global fields in the second one.

In the classical Langlands program, only the first two types of local and global fields are considered. The program, in its  simplest (and most studied) variant, consist in the  construction of  a correspondence between finite dimensional (of dimension $n$) representations of the Galois group $G_K = \Gal(K^{sep}/K)$ of a separable closure of the field $K$ and  irreducible infinite-dimensional (usually) representations of the group  $GL(n, A_K)$, where $A_K = K$ for the local fields $K$ and  $A_K = \A_K$ (ring of adeles) for the global fields $K$.

Recall that $\A_K = \{(f_v): \mbox{for all } v ~f_v \in K_v,~\mbox{for almost all } v ~f_v \in \hat{\cal O}_v  \}$. Here $v$ runs through all classes of valuations of the field $K$ with complete valuation rings $\hat{\cal O}_v$ and local fields $K_v = \Frac(\hat{\cal O}_v) \supset K$ which are the local completions of $K$. In the number field case, one  has to add finitely many archimedean valuations of the field  $K$ and the related embeddings of $K$ into the fields  $\R$ or $\C$.

Globally, the representations  that we consider here must be automorphic, i.e. realized in a space of functions on the quotient space
$\GL(n, \A_K)/\GL(n, K)$.  Since $n$-dimensional representations of the Galois group $G_K$ can be viewed as homomorphisms from $G_K$ in  $\GL(n, \C)$, then there is a natural generalization of this correspondence, when the group $\GL(n)$ is replaced by an arbitrary reductive algebraic group $G$.  In this more general situation, the hypothetical Langlands correspondence should look like (in a first approximation) as follows
$$
\begin{array}{llc}
Hom(G_K, ^LG(\C)) & \Leftrightarrow &
\{\mbox{irreducible  (automorphic)}\\ & & \mbox{representations of the group }~ G(A_K)\}
\end{array}
$$
Here  $^LG$ is a reductive algebraic group,  Langlands dual group to $G$ (see the definition and basic properties in \cite{T2, C}).
 It is important and highly nontrivial  that if the original group $G$ is regarded as a reductive group over the different  fields in the basic list, then the dual group  $^LG$ is always a group defined over the complex numbers\footnote{Or another field such as the field $\bar\Q_l$ which has nothing common with the ground field $K$.}. The group  $^LG$ contains the connected component of the identity $^LG_0$, which is a reductive algebraic group over $\C$.  Its root system is constructed by application of   duality of  tori, starting from the root system of the original group $G$. In addition, for the general linear group $G = \GL(n)$ the root systems of $G$  and  $^LG$  coincide. Thus, in this case, the presence of the duality is masked. Finally, the whole group $^LG$ is a semidirect product of groups  $^LG_0$ and the Galois group $\Gal(K^{sep}/K)$.

The representations of Galois groups included in the LHS of the Langlands correspondence require   substantial generalization, so that one can hope for the existence of bijections with automorphic representations of reductive groups that are on the RHS
 \footnote{We note that existence of the bijection was only assumed for the group  $G = \GL(n)$. Generally,  the correspondence may be not one-to-one and can have certain finite sets as preimages, i.e. so called    $L$-packets.}.
 In the theory of Artin $L$-series there are considered representations of finite Galois groups with values in the group $\GL (n, \C)$. If one goes  to the Galois group $G_K$ of a separable closure  of the  ground field $K$, then all of its continuous representations in the group $\GL (n, \C)$ are obtained from representations of finite factors  $\Gal (L/K)$, where $L/K$  runs through all the normal finite extensions of $K$. Thus, there are not too many  continuous complex representations of the Galois groups. On the other hand, there are continuous $l$-adic representations arising from Galois action on
etale cohomology of algebraic varieties defined over $K$. They usually do not pass through the finite factors of Galois groups \footnote{This is so if all cycles of the variety are algebraic. By the  Tate's conjectures,  the converse should be  true.}.
Deligne has invented how to build an extension of the Galois group $G_K$, whose complex representations  "coincide" with the $l$-adic representations of $G_K$. The Weil-Deligne group $WD_K$  is a semidirect product of the  Weil group $W_K$\footnote{The Weil group $ W_K $  of a local field $ K $ with finite residue field is the  subgroup of the Galois group $ G_K $, consisting of those elements which act on the residue field through  {\em integer} powers of the Frobenius automorphism.}
 with  $ \C$ and the action
$wzw^{-1} = q^{-\nu(w)}z$. Just,
$$
WD_K = \{w, z : w \in W_K,~z \in \C \},
$$
and
$$
(w, z)\cdot (w', z')  = (ww', z + q^{-\nu(w)}z'),
$$
where $\nu : W_K \rightarrow \Z$ is  the canonical mapping of the Weil group that sends $n$-th degree of the  Frobenius automorphism onto $n$
and $q$ is  the number of elements in the residue field \footnote{See details in \cite{T2}[n 4.1]. Instead of the Weil-Deligne group one can  consider  representations of the Weil-Arthur group $WA_K = W_K \times \SL(2, \C)$. The transition from representations of one group to the representations of another one  provided by the theorem of Jacobson-Morozov (see \cite{L3}).}.
\vskip0,5cm
In particular, the Langlands correspondence $\text{LC}_n$ for local fields and the group $\GL(n)$ must have the form:
$$
\begin{array}{llc}
Hom^{c, ss}(WD_K, \GL(n, \C)) & \Leftrightarrow &
\{\mbox{smooth irreducible} \\ && \mbox{representations of the group}~ \GL(n,K)\}
\end{array}
$$
and this is the final exact form (the representations of the Weil-Deligne group  supposed to be continous and completely reducible).
{\em Smoothness} of a representation $\pi: \GL(n, K) \rightarrow \text{End}(V) $ means that for every $ v \in V $ there is an open compact subgroup  $ \K'$ in $\GL (n, K) $ such  that $ \pi (\K')(v) = v$. {\em Irreducibility} means  the absence of {\em any} non-trivial invariant subspaces of $V$.
Let $ \rho, \rho '$ run over $n$-dimensional representations of the
Weil-Deligne group, and $\pi $ be a representation of $ \GL(n, K)$.  The Langlands correspondence has the following properties:
\vskip 1mm
\begin{itemize}
\item[i)] $\text{LC}_1(\Det\rho) = \text{central character of the representation } \text{LC}_n(\rho)$;
\item[ii)]$\text{LC}_n(\rho\otimes\chi) = \text{LC}_n(\rho)\otimes(\text{LC}_1(\chi)\circ\Det)$, where $\chi$ is a one-dimensional representation and $\Det : \GL(n, K) \rightarrow K^*$ is the determinant;
\item[iii)]$\text{LC}_{m+n}(\rho\oplus\rho') = \text{Ind}_P^{\GL(m+n, K)}\circ\text{Res}^P_M[\text{LC}_{m}(\rho) \boxtimes \text{LC}_{n}(\rho')]$;
\item[iv)]if the representation  $\rho$ is irreducible then representation  $\text{LC}_n(\rho)$ is a cuspidal one;
\item[v)]$\text{LC}_n(\check{\rho}) = \check{\text{LC}_n(\rho)}$;
\item[vi)] $L(\rho, s) = L(\text{LC}_n(\rho), s)$;
\item[vii)] if the representation $\rho$ is unramified then representation $\text{LC}_n(\rho)$ is a spherical one.
\end{itemize}
\vskip 1mm
For the property iii),  we have used the  embedding of the product  of groups $\GL(m, K) $ and $ \GL(n, K) $ into the group $ \GL(m + n, K) $ as the  Levi quotient $ M = {\GL(m, K) \times \GL(n, K)} $ of the standard parabolic subgroup $ P $. Concerning the operation $ \boxtimes $ see note. 18. The subsequent operation is {\em parabolic induction} which  extends  the representation of  the Levi subgroup $M$ to the  subgroup $ P $ and then induces it to  the whole group $ \GL(m + n, K) $ (for details, as well as for the definition of cuspidal representations, see \cite{BZ, C} and \cite{B} [ch. 4] for the case $ \GL(2) $).

Furthermore, $ \check\rho~ (\check\pi) $ denote  dual (contragradient) representations arising from the representations  $ \rho~ (\pi) $. Let  $ L(\rho, s) $ be the  Artin  $L$-function of  representation $ \rho $ (see \cite{T2}[n 3], \cite{B}[ch. 1.8]) and $ L(\pi, s) $ be  $ L$-function of representation $ \pi $ (see section on the Hasse-Weil conjecture ).  A representation  $ \pi $ of $ \GL(n, K) $ on a space $ V$  is called {\em spherical}, if there exists a vector $ v \in V$, such that $ \pi(\K)(v) = v $, where $ \K = \GL(n, {\cal O}_K) $ is a maximal compact subgroup. The list of properties can be further  supplemented by an equality for local constants appearing in the functional equation of  $L$-functions (see \cite{H}).

Similar properties should hold for the global Langlands correspondence.
\vskip 0,5cm
Up to now, the Langlands correspondence was constructed exactly in this form for the local fields of  type 1) (M. Harris-R. Taylor \cite{HT}; G. Henniart \cite{H}) and type 2) (G. Laumon-M. Rappoport-U.Stuhler \cite{LRS}). For the global fields of type  5) and the group  $G = \GL(n)$, the correspondence was obtained for $n = 2$ by V. G. Drinfeld \cite{Dr1, Dr2} and for arbitrary $n$  by L. Lafforgue \cite{L}.
In any case,  construction of the correspondence  in the number field  case remains an open question, and apparently very far from being resolved.

Previous constructions have a purely arithmetic origin. Based on his proof, which is concentrated around the moduli space
of vector bundles on algebraic curves, Drinfeld has expanded the class of fields and introduced the geometric Langlands
correspondence, which should take place for the fields of types 3) and 6).

The initial object will then be a complete algebraic curve $X$, defined over $\C$, and for which the field $K$ is a field of
rational functions on $X$. The correspondence looks as follows
$$
\mbox{Local}~ ^LG-\mbox{systems on}~ X~   \Leftrightarrow  \mbox{Hecke sheaves on scheme (stack)}~ Bun_{X, G}
$$
and was proven for unramified local systems and the group $G = \GL (n)$ (V. G. Drinfeld for $n = 2$ and E. Frenkel - D. Gaitsgory -
K.  Vilonen for any $n$, based on a construction of  G. Laumon). Since local (unramified) $^ LG$ system on $X$ are nothing but a representation
of the fundamental group of  $X$ in the group $^ LG$, it looks as a quite natural analog of the arithmetical situation.
On the other hand, there is $Bun_ {X, G}$, the  stack of principal $G$-bundles on $X$, and the ring of Hecke correspondences
acts on it. For the stack  $Bun_ {X, G}$, one needs not  the usual local systems but perverse sheaves which are considered instead of automorphic forms in the classical Langlands correspondence. This sheaf  is called the Hecke
sheaf, if it is semi-invariant with respect to Hecke correspondences. In this situation, the Hecke  sheaves do not
replace the representations of the adelic group responsible, according to the arithmetic Langlands correspondence, to the
representations of the Galois group, but they correspond to the spherical vectors of these representations. Whether
there are the representations itself, I do not know, but, if they exist, they should be in a reasonable sense unramified
(spherical)). In the arithmetical situation the spherical vector is uniquely defined, up to a constant, and generates
the representation space.
\vskip0,5cm
What are the problems in number theory that can be studied with help of  the Langlands program?

1. \underline{Non-abelian class field theory}.
 Langlands correspondence arose from consideration of the abelian situation
where $G = \GL (1)$ and the representation of the Galois groups are homomorphisms of $G_K$ into $\GL (1, \C) = \C ^ *$,
i.e. characters of $G_K$. Previously, such maps are called quasi-characters, keeping the term "character" for
the homomorphisms in the unitary group $U (1) \subset  \GL (1, \C)$. By class field theory we have the  reciprocity  map,
$$
C_K \rightarrow G^{ab}_K,
$$
where $C_K = K^ *$  for local fields of types 1) and 2) and $C_K = \A_K/K^ *$  for global fields of types  4) and 5).
Characters of $C_K$ (and they are automorphic forms on $\GL (1)$) can be transferred to the group $G_K$  as its abelian characters,
i.e. 1-dimensional representations. We thus obtain an abelian Langlands  correspondence \cite{L2}, which served as a sample and the starting point
of the whole program.

However, there is a fundamental difference between  the general Langlands construction and its abelian  case. Actually, the difference  is related to the problem of explicit description of  the Galois
group. Class field theory allows to compute explicitly the Galois group of the maximal abelian extension of the original
field (specifying its generators and relations). Moreover, the cohomological formulation of class field theory makes
it possible to compute the Galois group of some non-abelian extensions  ($l$-extensions \cite{Sh, Koch}).

The Langlands program  allows us to
describe the set of representations of the Galois group, providing it with a standard set of operations (primarily direct
and tensor product).  Further development of the program  was to construct a category of representations of
the Galois group as a monoidal category (more precisely, the Tannaka category or, probably,  its generalization, see discussion in \cite{L3, Ram}), which allows to reconstruct uniquely the group (in view of a general theorem
of the theory of Tannaka categories \cite{Br}). Note that this theorem is an existence theorem and it is unclear how to get an explicit
presentation through generators  and relations.
\vskip0,2cm
2. \underline{Conjectures on zeta- and $L$-functions} of arithmetical schemes (see below the section on the Hasse-Weil conjecture).
In this direction, a great achievement is the proof of Artin's conjecture on the holomorphy of $ L $-series for a wide class
of  two-dimensional representations of the Galois group over $ \Q $ (R. P. Langlands-J. Tunnell). Actually, a special case of the octahedral representations was  a starting point for  the A. Wiles proof of the Taniyama-Weil conjecture.
\vskip0,2cm
3. \underline{Description of}  the set of smooth irreducible  \underline{representations}\\
 of \underline{reductive groups } over local and global fields.
\vskip0,2cm
The last line is the most successful. We have seen how the point
of view of representation theory of reductive groups {\bf requires} certain changes of the Galois group on the Weil-Deligne group or the Weil-Arthur group.
This is most  sharply  shown by the example of a local field $K = \C$. The group $\GL (n, K)$ has for $K = \C$ a nontrivial set
of irreducible representations (all of them are the principal series representations and their  parameter space is a complex manifold = quotient-space of  the manifold of characters of  torus $(\C^*)^n$   by the permutation group $S_n$\footnote{D. P. Zhelobenko-M. A. Naimark's theorem (1966, see exposition in \cite{Zh}). Explicit form for the characters of the torus  $\C^*$ and consequently torus  $(\C^*)^n$ see below and more detailed in \cite{Kn}.}).  At the same time
Galois group $G_K$ is  trivial, but the Weil group  is, according to  the general theory \cite{AT}[ch. XV],\cite{T2}[n 1],   $W_{\C} = \C^*$ and has exactly the same space
of $n$-dimensional semisimple representations \cite{Kn}.

We now consider the basic fields of the Langlands  program in terms of general principles of arithmetic algebraic geometry.
\section{Basic fields from the viewpoint of the scheme theory}
%\begin{center}
%\vskip2,5cm
%{\large {\bf Basic fields from the viewpoint of the scheme theory}}
%\vskip0,5cm
%\end{center}
The fields, considered in the  Langlands' theory, are fields of functions on some schemes. In the scheme theory, we have a
classification by  a dimension. If the scheme is affine
and corresponds to a ring $A$, then we  can take  its Krull  dimension. We have
the following table of rings, arising  in arithmetic:
$$
\begin{tabular}{|c|c|c|} \hlx{hv}
$\Dim(A)$ & geometric case & number case \\ \hlx{vhv}
$>2$ & $\ldots$ & $\ldots$ \\ \hlx{vhv}
$2$ & $\F_q[x,y]$ & $\Z[y]$ \\ \hlx{vv}
& (surfaces) & (arithmetic surfaces) \\ \hlx{vhv}
$1$ & $\F_q[x]$ & $\Z$ \\ \hlx{vv}
& (curves) & (arithmetic curves, $\Spec(\Z)$) \\ \hlx{vhv}
$0$ & $\F_q$ & $\F_1$ \\ \hlx{vh}
\end{tabular}
$$
Here, $\F_q$ is  a finite field of $q$ elements, $\Z$ is  the ring of integers
and $ ~ \F_1$ is  now popular field of one element.

These rings correspond to affine schemes and the fields of functions on them.
In particular, this table
includes global fields of types 4) and 5) from the list of fields that are now considered
as  the main ones for the Langlands program (see  ~ \cite {L3, F2}):
$$
\begin{tabular}{|c|c|c|} \hlx{hv}
$\Dim(A)$ & geometric case & number case \\ \hlx{vhv}
$2$ & $\F_q(x,y)$ & $\Q(y)$ \\ \hlx{vv}
& (rational functions) & (rational functions) \\ \hlx{vhv}
$1$ & $\F_q(x)$ & $\Q$ \\ \hlx{vv}
& (rational functions) & (rational numbers) \\ \hlx{vhv}
$0$ & $\F_q$ & $\F_1$ \\ \hlx{vh}
\end{tabular}
$$
Finally, there is a table of the local fields that arise under these circumstances:
$$
\begin{tabular}{|c|c|c|} \hlx{hv}
$\Dim(A)$ & geometric case & number case \\ \hlx{vhv}
$2$ & $\F_q((x))((y))$ & $\Q_p ((y))$, \, $\R((y))$, \, $\C((y))$ \\ \hlx{vv}
& (iterated Laurent power series) & (Laurent power series) \\ \hlx{vhv}
$1$ & $\F_q((x))$ & $\Q_p$, \, $\R$, \, $\C$ \\ \hlx{vv}
& (Laurent power series) & ($p$-adique, real and complex \\ \hlx{vv}
&& numbers) \\ \hlx{vh}
\end{tabular}
$$

Here we see all types  1) - 3) of local fields, but also a number of the new fields which were introduced in the higher adelic theory \cite{P3}\footnote{See their informal definition at the beginning of the following section.}.
Such a classification is the development of  ideas of Andre Weil from his  report to the Cambridge congress in 1950 ~\cite{W2},  where he stressed the importance of the concept of the Krull dimension
for the classification of arithmetic problems. In his pre-war letter to  sister
Simone ~\cite{W1} where Weil also speaks on the analogy between the number
and the geometric case, the situation is simplified, and so he adds the field
$\C((y))$  and puts it on the "one line"  with such fields as $\F_q ((x))$,
$\Q_p$, what as seen from the table is not consistent with the dimension.

Returning to the global fields, we may ask, where in these tables are
such schemes as  the curves over the complex numbers associated with fields of type 6)? These curves and the fields of
functions on them
are the subject of the geometric Langlands correspondence, generalizing the usual
Langlands correspondence.
Indeed, they are  not here, as they correspond to the
rings of the form $\C [x]$. These objects have an intermediate nature, they are global
and of dimension $1$ over $x$ (the coordinate on the curve), and local in the field of definition  $\C$.
This field is included in the table as a one-dimensional (!) local field. Thus, such a ring as $\C [x]$
correspond to a geometric object of dimension $2$. That is, they are archimedean  fibers of arithmetic surfaces
(on them  see below).

This observation and the fundamental
arithmetic analogies \cite{P2}  suggest to consider, as a "partner" of
the geometric Langlands correspondence for curves over $\C$,
constructions of this kind for curves over $p$-adic numbers $\Q_p$. As far as I know,
this was not considered by anyone and later  we'll discuss how it might look\footnote{As was pointed out to me by \.Ilhan \.Ikeda, a construction of this kind was recently found by A. G. M. Paulin for abelian case.}.

Summarizing this discussion, we can say that the list of six basic fields could be completed as follows. Only four types of fields - 1), 2), 4), 5) belong to the  one-dimensional situation. Fields of types 3) and 6) are already in a two-dimensional situation. In this case we have two types of local fields
\vskip2mm
1) finite field extensions of $ {\Q} _p ((t))$, $ {\R}((t)) $, $ {\C}((t)) $ and $ {\Q} _p \{\{t\}\} $ \footnote {If $K$ is locally compact field with a non-Archimedean norm $ \vert \cdot \vert $, then the field $ K\{\{t\}\}$ consists of the infinite in both directions series $\sum_i a_i, ~ a_i \in K $ with
$ \vert a_i \vert \leq O(1), ~ \vert a_i \vert \to 0 $ for $ i \to - \infty $. The  fields $ {\Q} _p \{\{t\}\} $ can arise on the arithmetic surfaces.}.

2) the field of Laurent power series $ {\F} _q ((u))((t)) $.
\vskip2mm
and two types of global fields:
\vskip2mm
3) the fields of functions on arithmetic surfaces (= finite extensions of the field $\Q(x)$ of rational functions).

4) the fields of functions on algebraic surfaces over a finite field of constants (= finite extensions of the fields $\F_q (x, y)$ of rational functions in two variables).
\vskip2mm
Between them lie the fields of intermediate types (partially  local, partially global). In the situation of a surface over a finite field the  intermediate fields are
\vskip2mm
2a) finite  extensions of the fields ${\F} _q (u)((t))$.

2b)  the field of fractions of the ring ${\F} _q [[u, t]]$ of formal power series in two variables.
\vskip2mm
In the case of arithmetic surfaces, the field $\C(x)$ (the field 6) from the list of Langlands) is also an intermediate field
but it is not included in the standard set of fields appearing in the two dimensional adelic theory (see below and \cite {P3}). The reason is that  it due to the additional structure: the existence of a morphism from the two-dimensional scheme onto a one-dimensional scheme  $B$ and its fibers, defined over  one-dimensional local  fields of the scheme $B$.

{\bf Remark 1}. Adelic technique allows us to consider local and global fields of the original Langlands program  in a simultaneous
way. As we have seen, in both cases, these fields are the fields of  functions on a scheme $X$. In the local (non-Archimedean!) case
we have $X = Spec({\cal O})$,  where ${\cal O}$  is a local discrete valuation ring with finite residue field. In
global situation,  $X$ is either a curve over a finite field or spectrum of the ring of an algebraic number field.
In both cases,  $A_K = \A_X$  if we apply the general definition of the ring of adeles \cite{FP,P3},  which may  include and the archimedean components as well.
\section{Two-dimensional generalization of the Langlands correspondence}
%\begin{center}
%\vskip 2,5cm
%{\large {\bf Two-dimensional generalization \\ of the Langlands correspondence}}
%\vskip 0,5cm
%\end{center}
We may wander about  existence of  a two-dimensional generalization of the Langlands correspondence.
Let $X$ be a two-dimensional scheme  and $K$ be the field of functions on it. Then the  $n$-dimensional representations
 of the Galois group $G_K$ of a separable closure for  $K$ could
correspond to irreducible automorphic representations of the adelic group $G(\A_X)$.
Here, the ring of adeles $\A_X$ was introduced in the higher adelic theory in the 1970s
(see ~ surveys  \cite{FP, P3}).

If $ x \in C \subset X $ is a flag on $ X $, consisting of a point $ x $ and an irreducible curve $ C $, and the point $x$ is smooth on $C$ then one can introduce
 a two-dimensional local field $ K_ {x, C} $. In local coordinates $ u, t $ of a formal neighborhood
 of the point $ x $, the  field $ K_ {x, C} $ is equal to $ k (x) ((u)) ((t)) $ (if $X$ is a smooth surface),   where  $ C = (t = 0) $ and  $ k (x) $ is the residue field of
 $ x $. In general, one can attach a direct sum of finitely many two-dimensional local fields to the flag $x, C$. The adelic ring
 $ \A_X $ is defined as the adelic (a part of the total)  product of  the fields $ K_ {x, C} $ for all flags $ x, C $.

The field $ K_ {x, C} $ contains a discrete valuation ring $ \hat{\cal O}_ {x, C} $ (which is $ k (x) ((u)) [[t]] $ if $X$ is a surface).  For an irreducible curve
$C$, we take $K_C = \Frac(\hat{\cal O}_C)$ and for a point $x$, we set $K_x =  K\hat{\cal O}_x \subset \Frac(\hat{\cal O}_x)$.
Here, $\hat{\cal O}_C = k(C)((t))$ and $\hat{\cal O}_x = k(x)[[u, t]]$.
The adelic product of these  rings give rise to three subrings $\A_{12}, \A_{01}, \A_{02}$ in $\A_X$. These rings are embedded into the $\A_X$ in a diagonal way (exactly as the principal adeles in the case of one-dimensional schemes).
 We can compare the structure
of the local adelic components in dimensions one and two:
$$
\xymatrix@!0{
K_x\ar@{-}[dd] & & &  & K_{x,C}\ar@{-}[dl] \ar@{-}[dr]  &
\\
 & & & K_x\ar@{-}[dr]  & & K_C \ar@{-}[dl]
\\
K & & & & K &
}
$$

A generalization of the Langlands program that uses this
theory of adeles, were outlined by  M. ~ Kapranov \cite{K} (and detailed
in ~\cite{I}). It uses not ordinary notion of representation, but  its
generalization  associated with 2-categories \cite{GK}.

Here the number 2 stands as  dimension of the scheme $ X $. We have the following table of categories:
$$
\begin{tabular}{|c|c|c|c|} \hlx{hv}
$n = \text{dim}$ & $n$-categories & objects & distinguished  \\ \hlx{vv}
&  & & object \\ \hlx{vhv}
$2$ & $\text{2-Vect}/k$ & $\text{Vect-mod}/k$ &  $\text{Vect}/k$ \\ \hlx{vhv}
$1$ & $\text{Vect}/k$ & $V/k$ &  $k$\\ \hlx{vhv}
$0$ & $k$ & $A \subset k$ & $\{1\}$ \\ \hlx{vhv}
\end{tabular}
$$
where   let $k$ be a ground field (or a ring), $A$ be its subset,  $V$ be vector spaces over $k$ and  $\text{Vect}/k$ be the category of all vector spaces over  $k$.
The latter is a tensor category, and it can be considered as a categorical generalization of the notion of  ring. Then, we can
 define  a category of modules $ \text{Vect-mod}/k $ over  such a "ring" and they  form a 2-category whose objects are
  such categories.

We see an obvious inductive structure that can be easily compared with   inductive structure of the $n$-dimensional local fields.
 In addition, if there is a group $ G $, then we can introduce a concept of $ n $-representation. 1-representation  is the usual
 homomorphism in the group of automorphisms of a vector space. 2-representations are obtained  replacing the spaces by  categories
 and the homomorphisms by   functors. More precisely, instead of the space $ V $ we choose the category module $\text{Vect-mod}/k $ and the 2-representation $ \pi $ is to attach to an element $ g $ of $ G $ a  functor $\pi (g) :  \text{Vect-mod}/k \rightarrow \text{Vect-mod}/k$,
 satisfying a condition of multiplicativity. This generalization is a special example of the general categorification process of
 concepts and structures\footnote{ The 0-representations will be discussed later in appendix.}.

Then, in a certain approximation the Langlands correspondence for two-dimensional local field $ K $ looks, according  Kapranov,  as follows:
$$
\begin{array}{ccc}
\{n\mbox{-dimensional complex (or } l\mbox{-adic)} &  \Leftrightarrow  &  \{\mbox{irreducible  2-representations}\\
\mbox{Galois representations of }~ G_K\} && \mbox{of the group  }~ \GL(2n,K)\}
\end{array}
$$
 and for a global field $ K $ on the scheme $ X $ the field $ K $ in the RHS side of the correspondence is replaced by the ring $ \A_X $ of adeles.
 The  2-representation must satisfy a certain automorphic  condition. These definitions and related concepts need to be clarified and one has to understand certain restrictions they hold.
 In particular,
  the characteristic of the basic field is supposed to be 0 and nothing suggested for the local fields of archimedean type.
%  Also, the representation of the Galois group we have to consider not just any, but with a condition on the simple
%components such as the multiplicity 1 of the  spectrum.
   Surely, the construction has to exist  for these fields and for the fields of finite characteristic. We shall return to this issue \footnote{Recently, D. V.  Osipov, developing the ideas of Kapranov, proposed a definition of unramified local Langlands correspondence for two-dimensional local fields and the group $ \GL(n) $ for any $n$ \cite{O2}.}.
 \section{Functorial properties    of the Langlands correspondence}
% \vskip1cm
% \begin{center}
%\vskip 0,5cm
%{\large {\bf Functorial properties   \\ of the Langlands correspondence}}
%\vskip 0,5cm
%\end{center}
The usual Langlands correspondence has a large number of functorial
 properties. The theory includes two kinds of objects: a scheme or a field and a reductive algebraic group (better to say, a root data).
The Langlands principle of functoriality \cite{L3} says what happens when one changes the group. We will restrict ourselves here with the properties related to a change of the schemes\footnote{These two are not independent. In fact, the functoriality principle implies the automorphic induction and base change for finite extensions.}. First, they are the properties associated with finite extensions
of fields (base change and automorphic induction or lifting).
Second, the properties associated with the transition from local to global fields.
All these properties are special cases of the general functorial construction of the inverse
and direct images.

The first properties  related to  finite surjective morphisms of one-dimensional schemes
$f\colon X \to Y$, where both $X, Y$ are either curves (or rings of integers in algebraic number fields)
or the spectrum of local rings.
The latter related to  morphisms $f\colon X \to Y$, where $X$  is a spectrum of a  local
ring, and $Y$ is a  curve.
Let $B$ be a curve over $\F_q$ and $b \in B$. Then there is the commutative diagram
$$
\xymatrix@C+1pt{
\text{sheaves ${\cal F}$  on  } B \ar[r] \ar[d]_{i^*} & \text{ automorphic forms on } \GL(\A_B)
\ar[d]_{i^*} \\
\text{sheaves } i^*({\cal F}) \text{ on } \Spec({\hat {\cal O}}_b) \ar[r] &
\text{ smooth functions on } \GL(K_b),
}
$$
where $i : \Spec({\hat {\cal O}}_b) \rightarrow B$ is the canonical embedding, the left map $i^*$ is the inverse image of sheaves and the right
map $i^*$ is  taking of the local $b$-component of a function on the whole adelic group (the non-trivial examples of this diagram see in \cite{HT}).

If the  Langlands  correspondence exists in the one form or another for the higher
dimensions, it should have a much greater number  of
functorialities. Indeed, let a two-dimensional scheme $X$ be considered together with
a structure morphism $f\colon X \to B$, where $B$ is an one-dimensional scheme.

Let us assume for simplicity that our schemes are varieties defined over a finite field. For $l$-adic sheaves (representations of the Galois groups), there are direct (pushforward) and inverse (pullback)
images under the mapping $f$. Then, by Langlands correspondence, the same kind of operations should
exist   for representations of adelic  groups.

Without discussing this question in full generality, let us consider a special and already extremely
interesting case. Let us start with an  one-dimensional representation of the Galois group on  $X$.
It determines certain  sheaf  ${\cal F}$ of rank 1 on $X$ and the direct images $R^i f_*({\cal F})$
on  $B,~i = 0, 1, 2$.

In this abelian situation on the surface $X$,   there is already the Langlands correspondence between
characters of the Galois group $G_K$ and characters of  the group $K_2 (\A_X)$. By the
two-dimensional class field theory (see surveys \cite{FP, Inv, R}), there is a  canonical map
$$
\varphi_X :  K_2(\A_X) \longrightarrow G_K^{ab},\eqno (1)
$$
which gives the map of groups of characters in the opposite direction.
The reciprocity laws on the $X$ read as
$$
\varphi_X (K_2(\A_{01})) = (1), ~  \varphi_X (K_2(\A_{02}))  = (1). \eqno (2)
$$
We can assume that on the "automorphic side" of the Langlands correspondence
 there are also such operations as direct  and inverse images  between the characters of $K_2 (\A_X)$
and automorphic forms on $\A_B$.

The case of the {\bf inverse image} is not hypothetical, but follows from the known
results of the higher adelic theory. There is an operation of the direct image for abelian
groups~\cite{Ka, O1}\footnote{
In this paper, the map $f$ is smooth. The analogous homomorpism for differential forms was constructed in \cite{M} for non-smooth maps.
 Quite  recently \cite{Liu},  the direct image morphism of the $ K$-groups  was built  for non-smooth maps of algebraic and also arithmetic surfaces.
}.
$$
 f_* \colon K_2(\A_X) \to K_1(\A_B).  \eqno (3)
$$
Since $K_1 (\A_B) =
\GL(1, \A_B)$, we obtain a map of automorphic forms on $\GL(1, \A_B)$
to the  automorphic forms, that is, ~ characters on $K_2 (\A_X)$. Their "automorphy"
will be justified  by higher reciprocity laws   on the scheme ~ $X$.
We have the commutative diagram
$$
\xymatrix{
{\begin{tabular}{|c|}\hline Sheaves \\ $f^*({\cal F})$ \\  on  $X)$\\ \hline\end{tabular}}  \ar[r]
 & {\begin{tabular}{|c|}\hline Characters\\ of the group\\   $K_2(\A_X)$\\ \hline\end{tabular}}  \\
{\begin{tabular}{|c|}\hline Sheaves  ${\cal F}$ \\ of rank 1 \\ on $B$\\ \hline \end{tabular}}\ar[r] \ar[u]_{f^*}&
 {\begin{tabular}{|c|}\hline Automorphic \\ forms \\ on $\GL(1, \A_B)$.\\ \hline \end{tabular}}\ar[u]_{f^*}
} \eqno (4)
$$
The {\bf direct image} must\footnote{We deliberately describe this hypothetical construction in a preliminary vague form. Generally, the Langslands correspondance  includes the irreducible automorphic representations on the RHS. Nevertheless, in many cases it is possible to find a generating vector in the spaces of representations. These vectors lead to  automorphic forms on the adelic groups. So, for an  unramified case there  are spherical vectors and unramified automorphic functions; for generic representations one can use the Whittaker functions. At the end of this note we give an exact formulation in a simple but still non-trivial case. }
consist in a construction of automorphic forms on $\GL(r_i, \A_B)$,
corresponding to the characters  of  $K_2 (\A_X)$ on $X$. The following
commutative diagram should be true:
$$
\xymatrix{
{\begin{tabular}{|c|}\hline Sheaves ${\cal F}$\\ of rank 1  \\  on  $X$\\ \hline\end{tabular}}  \ar[r] \ar[d]_(.45){R^if_*}
& {\begin{tabular}{|c|}\hline Characters\\ of the group \\   $K_2(\A_X)$\\ \hline\end{tabular}}  \ar[d]_(.45){f^i_*} \\
{\begin{tabular}{|c|}\hline Sheaves \\ $R^if_*({\cal F})$\\ on $B$\\ \hline \end{tabular}}\ar[r]&
  {\begin{tabular}{|c|}\hline Automorphic \\ forms\\ on $\GL(r_i, \A_B)$.\\ \hline \end{tabular}}
} \eqno (5)
$$

Here $r_i$ is  rank of the sheaf $R^i f_* (F)$. It is quite possible that instead of the groups
$\GL(r_i, \A_B)$ separately, there must be considered "even" $~\GL (2, \A_B)
\sim \GL(r_0, \A_B) \times \GL(r_2, \A_B)$ and "odd" $~\GL (r_1, \A_B)$
groups. Also,  in the case of "middle" group, $i = 1$, the monodromy  preserving the skew-symmetric cup product could be
taken into account and possibly one has to consider  the orthogonal group corresponding to symplectic group (of the Galois action) according to the $L$-duality \footnote{More exactly, the group $\text{GSpin}(2g + 1, \C)$ which is $L$-dual to the group $\text{GSp}(2g)$ of simplectic similitudes.}.

{\bf Remark 2}. We may compare this picture with a number field case
 where a two-dimensional regular  scheme $X$
is mapped on a one-dimensional scheme $B$ (= the ring of integers in a number field). We expect the same kind of construction also in this case.
  If the general fiber of the mapping $ f $ is an elliptic curve over $ \Q $, $B = \Spec(\Z)$ and $ \chi $ is  trivial,
 then the  automorphic form $ f^1_*(\chi) $ is nothing but an  adelic version of the parabolic form corresponding to  the curve
 according to the  Taniyama-Weil conjecture (now a theorem!). In the  Langlands program, this form generates  a cuspidal representation
 of the group $ \GL(2, A_B)$ \cite{L1, G}. At the same time, the  less complicated  representation associated
 (in our notation) with  forms  $f^0_*(\chi)$ and $f^2_*(\chi) $ somehow remained in the shadows. These forms are defined on the group $ \GL(1, \A_B) $.
Thus, they are  characters of $ \A^*_B $ and, therefore, they give a character $ \eta $ of the standard  maximal torus $ T $ in $ G = \GL(2, \A_B) $. Parabolic induction applied to  $ \eta $ (if it is  not equal to 0) gives  an  induced representation from  the principal series of $ G $ which is  generated by
an  Eisenstein series. Thus,  the original character $ \chi $ can be associated with {\bf  two} infinite-dimensional representations
 of  $ \GL(2, \A_B) $,  the "even" one belonging to the continuous spectrum, and the other "odd" one belonging to the discrete spectrum.
\vskip 3mm
Return back to our diagram.
Since we consider  the case of curves over a finite field the  correspondence of  the bottom line
was constructed (Drinfeld-Lafforgue) and one "only" needs  to close the entire diagram.
It seems natural to build  first the local and semi-local direct
images that enter into  the following  diagram
$$
\xymatrix{
{\begin{tabular}{|c|}\hline Characters\\ of the group\\   $K_2(\A_X)$\\ \hline\end{tabular}}  \ar[r] \ar[d]_(.45){f^i_*} & {\begin{tabular}{|c|}\hline Characters\\ of the group\\   $K_2(\A_{X_b})$\\ \hline\end{tabular}}
 \ar[r] \ar[d]_(.45){f^i_*} & {\begin{tabular}{|c|}\hline Characters\\ of the group\\   $K_2(K_{x,C})$\\ \hline\end{tabular}}  \ar[d]_(.45){f^i_*} \\
{\begin{tabular}{|c|}\hline Automorphic\\ forms\\ on $\GL(r_i,\A_B)$\\ \hline \end{tabular}}\ar[r]&
 {\begin{tabular}{|c|}\hline Automorphic\\ forms\\ on $\GL(r_i,K_b)$\\ \hline \end{tabular}}\ar[r]&
 {\begin{tabular}{|c|}\hline Automorphic\\ forms \\ on $\GL(r_i,K_b)$.\\ \hline \end{tabular}}
} \eqno (6)
$$
%$$
%\xymatrix@C-17pt{
%\text{Характеры } K_2(K_x,C) \ar[d]_{f^i_*} & \text{ Характеры } K_2(\A_{X_b})
%\ar[l] \ar[d]_{f^i_*} & \text{ Характеры } K_2(\A_X) \ar[l] \ar[d]_{f^i_*} \\
%\text{ Автоморфные формы } & \text{ Автоморфные формы } & \text{ Автоморфные
%формы } \\ \vspace{-20mm} \text{ на } \GL(m_i,K_b) & \text{ на } \GL(m_i,K_b) &
%\text{ на } \GL(m_i,\A_B).
%}
%$$
Here, $b = f (x)$, $C$  is a curve on $X$, passing through $x$ (in particular, the  fiber $F_b$ of the map $f$), $X_b$ is the
two-dimensional scheme, $X_b = X\times_B \Spec ({\mathcal O}_b)=$  infinitesimal
neighborhood of the fiber $f^{-1}(b)$.

The next step should be to prove the global automorphic property for the
 image-form on  $B$, starting from the  reciprocity laws for the character to $K_2 (\A_X)$.
Note that these issues can be considered completely independent of how one can  arrange
the Langlands correspondence for sheaves ${\cal F}$ of  rank > 1 on $X$ .

Also, it is  important to consider  functorial properties of  the Langlands correspondence
 for one-dimensional schemes (conventional theory) and schemes (= finite fields) of
dimension 0. This question was raised by me in the late 1970s  and discussed
then with A. ~ Beilinson. Although it was clear that the correspondence in dimension 0
must be something simple,  we could not think it out. When M. ~ Kapranov
sent me the first version of his work, I suggested that he has to find  this definition  and add it to his    two-dimensional construction and he did this in the final text. However, he did not consider  functorial properties associated with  surjective maps of the schemes of {\bf different} dimensions.

In the case of schemes of dimension 1 and 0, we have a diagram
$$
\xymatrix@C+1pt{
\text{sheaves ${\cal F}$ of rank } r \text{ on } X \ar[r] \ar[d]_{R^i f_*} & \text{ automorphic forms on } \GL(r,\A_X) \ar[d]_{f^i_*} \\
\text{sheaves   } R^i f_*({\cal F}) \text{ on }
\Spec (\F_q) \ar[r] & \;L\text{-functions},
}
$$
where  $X$ is a curve over a finite field $\F_q$, sheaves on  $\Spec (\F_q)$ are vector spaces with Frobenius action and
 $L$-functions $=$
characteristic polynomials of the Frobenius automorphisms.
 The lower line refers to zero-dimensional
situation and is the Langlands correspondence for dimension zero according to
Kapranov‘s proposal\footnote{See a discussion of this proposal below in appendix.}.

This diagram should be consistent with the hypothetical diagram which we  proposed
above  for a morphism of scheme $X_b$ on $\Spec ({\cal O}_b)$, if we
take a closed fiber $F_b = X \times \Spec(k(b))$ as the curve $X$. The starting point of such concordance is
the cartesian commutative diagram
$$
\xymatrix@C+1pt{
F_b \ar[r] \ar[d] & X_b \ar[d] \\  Spec(k(b)) \ar[r] & \;Spec({\cal O}_b).
}
$$
\section{Relation with the geometric Drinfeld-Langlands correspondence}
%\begin{center}
%\vskip 0,5cm
%{\large {\bf Relation with the geometric Langlands correspondence}}
%\vskip 0,5cm
%\end{center}
Let us now  consider the same picture for the  number field case where the scheme $X$
is mapped onto an  one-dimensional scheme $B$ (= the ring of integers) and completed in accordance with
the Arakelov theory by  archimedean fibers $X \times_B \C$ or $X \times_B \R$ \cite{P2}.

Semi-local schemes $X_b$ in the new archimedean situation no longer exist, but we
{\bf must} find an analogue of the above structures and in this case as well (according to
the general principles of analogy "number $\sim$ functions" in arithmetic ~\cite{P2}).
It would be natural to expect that such a two-dimensional (!) Kapranov's analogue of  the Langlands correspondence for
scheme $X \times_B \C$ (so far only very hypothetical) is already formulated
and partially proven  geometric Drinfeld-Langlands correspondence on the curve $X$,
defined over the complex field ~ $\C$, and which appears here as an archimedean fiber.
We will discuss this possibility later.

Let us  assume now that in this situation, there are also direct images
and the commutative diagram
$$
\xymatrix@C+1pt{
\text{sheaves ${\cal F}$  of rank 1 on } X \ar[r] \ar[d]_{R^i f_*} & \text{ characters $\chi$ of } C_X
\ar[d]_{f^i_*} \\ \text{sheaves  }~ R^i f_*({\cal F}) \text{ on  } \Spec (\C) \ar[r] &
\text{ automorphic forms on  } \GL(r_i,\C)
}
$$
where $f \colon X \to \Spec (\C), C_X$ is a  suitable class formation  (see \cite{AT}) in the
class field theory for the  field $\C (X)$ of functions on $X$. Such formations are constructed together with
class field theory for the fields of type $\C(X)$ in the 1950s  by  J.-P. Serre,
M. Hasewinkel ~\cite{S2, Has} and Y. Kawada-J. Tate ~\cite{KT} . They are as follows,
$$
C_X = \pi_1(\A_X^{*(0)}/K^*) = \pi_1(\lim_{\longleftarrow} J_m(X)),
$$
where $(0)$ denotes the degree $0$ ideles,  $J_m(X)$  are  generalized Jacobian varieties with "module" $m$ \cite{S1}, the limit
is taken over all
$m$ and it is a pro-algebraic group in the sense of Serre \cite{S2}.
%$lim_{\to}_{i \in I}~ lim_{\gets}_{j \le i}$
\vskip 0,5cm
The bottom line of the diagram  is well known in the classical theory of the
local Langlands correspondence for the  field $\C$. "Sheaves" on $\Spec(\C )$ are
vector spaces of finite dimension over $\C$ equipped with a  Hodge  structure. Note that the cohomology of the
local system (= sheaf) on $X$ can be equipped
with a canonical Hodge structure. And it is this structure that defines a representation
of the Weil group $W_{\C}$ for the field $\C$, which plays a role, in this situation, of the Galois group.
We have $W_{\C}= R_{\C/\R}\G_m = \C^*$.
So, there   arises
\vskip 0,5cm
{\bf Question  1}. How to construct in an explicit way the direct image =  the RHS arrow of the diagram?
\vskip 0,5cm

But one  first need to figure out what is the correct definition of the Langlands correspondence in this situation.
As above, we can and must begin with the local Langlands correspondence, local on $X$.
The existing Drinfeld theory (see a survey \cite{F1}) consideres a local correspondence  associated with the local fields of the form $\C((t_P))$ and that
correspond to points $P$ of $X$ (Here, $t_P$) is a formal local parameter at the point $P$). The local correspondence  must describe the $n$-dimensional representations of the group,
which is a local version of the Galois group (or of the fundamental group in the unramified case).
This is  already  interesting to study the case of an abelian group and $n = 1$. Thus we arrive at  the question about
the local analogue of class field theory for such fields. Since they are archimedean variants of two-dimensional
 local fields, it is natural to ask:
\vskip 0,5cm
{\bf Question 2}. What is the relationship with the two-dimensional adelic class field theory?
\vskip 0,5cm

If $K$ is a two-dimensional local field with finite last residue field, then the class field theory consists
in a construction of the canonical reciprocity map,
$$
K_2(K) \longrightarrow G_K^{ab},
$$
which is a local component of the global map considered above. The kernel of this map coincides with the subgroup of  divisible elements in  $K_2(K)$ (I. Fesenko's theorem, see his survey in  \cite{Inv}).
Examples of two-dimensional local fields are the fields $\F_q((u))((t)),~\Q_p((t))$ and for them the class formations may be determined
(in a very non-trivial way) and the Weil groups are defined  as an extension  of the Galois group by a group  $K_2(K)$ \cite{Ko, A}.

In the case of local field $\C((t))$ with archimedean residue field,
we have the following structure of $K$-groups
$$
K_2\C((t)) \rightarrow K_1\C = \C^*,
$$
where the group $ \C^* $ is not the Galois group of $ \C $, but   its Weil group (!) and one may assume that the Weil group
of the field $\C((t))$ can be defined by the same construction as above, making use of  the group $K_2 \C((t))$. This gives us the LHS
part of the local abelian Langlands correspondence for  $\C((t))$. The question of the RHS side of correspondence may be linked (according to
Kapranov) to a study of 2-representations.

In any case, the resulting structure is basically different from the local geometric  Drinfeld correspondence, which includes
connections on  line bundles, $D$-modules  and the Fourier-Mukai transformation (see survey \cite{F1}).
The fundamental point of this difference is the fact that the proposed construction takes into account the arithmetic
 nature of the ground field $ \C $ which is manifested in  that it may be a completion of a  global number field in
 some infinite point. When a curve  over $\C$ is an archimedean fiber of an arithmetic surface $X$,  the field $\C$ can be identified with
 the local completions of  the number fields related to all horizontal curves on $Y$. The  picture in  \cite{P2}[section 4]  is a simplest example.
%\begin{figure}[t!]
%\centering
%\includegraphics{picture}
%\caption{}
%\label{fig1}
%\end{figure}
%\par\smallskip
%\includegraphics{3.eps}
%\par\smallskip
 There, on the arithmetical surface $X = {\mathbb P}^1/\Spec(\Z)$  we see three horizontal curves defined by the equations $t = 0; t = 1/5; t = 2$ where $t$ is a coordinate on the general fiber ${\mathbb P}^1/\R$ which is the curve lying over the infinite place $(\infty)$ of the ground "curve" (= "compactification" of  $\Spec(\Z)$).

 At the same time, the Drinfeld construction  makes sense  over any ground field (of characteristic 0 ?) and does not depend on
 its nature.

\vskip 2mm
{\bf Remark 3}.
Also, there should  be a similar picture  for curves over ~ $\R$. There arise the local fields of the form
$\R((t))$ and for them we have
$$
\begin{array}{l}
K_2\R[[t]] \hookrightarrow K_2\R((t)) \rightarrow K_1\R = \R^*,\\
K_2\R[[t]] \rightarrow K_2\R \rightarrow \{ \pm 1\}.
\end{array}
$$
These expansions correspond to unramified extension  $ \C((t)) \supset \R((t)) $ and to tamely ramified extension
 $ \R((t^{1/2})) \supset \R((t)) $. There are no other abelian extensions of the  field $ \R((t)) $. The class field theory for the field
 $K = \R((t))$ is a canonical  isomorphism
 $$
 K_2(K)/\{\text{divisible elements}\} \rightarrow \Gal(K^{ab}/K).
 $$
 As above, we assume that the Weil group of the field $\R((t))$ can be constructed from  the group $K_2 \R((t)))$.

Consequently, the basic  fields of the Langlands  program should be
{\bf supplemented}  by finite  extensions of the fields  $\R(t)$ and $\Q_p(t)$, associated with algebraic curves defined over
fields $\R$ and $\Q_p$.
To the list of local fields, we have to add fields $\R((t))$ and $K((t))$ with $K \supset \Q_p$.

\vskip 2mm

{\bf Remark 4}. In the Langlands correspondence, a  fundamental role played by $ L $-functions. One can attach to  the representations of  Galois groups (and more general of the  Weil groups) certain  $L$-functions (in particular, these are the Dedekind and Hecke $ L $-functions for  abelian representations and their
generalization to the non-abelian case, introduced by E. Artin). On the other hand, R. Langlands (and then
R. Godement and H. Jacquet) introduced  the $ L $-functions for  automorphic representations
of adelic groups of reductive groups. As the arithmetic $ L $-functions, they
determined by an Euler product of local $ L $-factors,
associated with the local components of a global automorphic representation (see the section on  Hasse-Weil conjecture).

An important property of  the Langlands correspondence is the equality of $ L $-functions which stand in both parts
of the correspondence, arithmetic $ L $-functions and automorphic $ L $-functions. Global equality precedes by equalities of local $ L$-functions  for all the local correspondences.

In the case of  geometric Langlands correspondence, introduced by Drinfeld,
apparently there is no proper definition of
$ L $-functions and hence there is no any  equality of this kind. Our interpretation of geometric correspondence suggests that there have to
exist $ L $-functions associated with the local fields $ \C((t_P)) $ of points $P$ on
the curve $ X $ defined over $ \C $. They should be included in the local
geometric Langlands correspondence. Since there is a structure map $ X \rightarrow \Spec(\C) $, one needs to take into account that the field $ \C $ has its own  $ L $-function entered  into  the local Langlands correspondence for this field. They  are the gamma functions, associated with Hodge structures according to Serre  \cite{S3}.

Since we expect that there are certain direct and inverse images of automorphic forms associated with the map $ f: X \rightarrow \Spec(\C) $, this  suggests that the hypothetical
 $ L $-functions of local fields $ \C ((t_P)) $ for the points on  $ X $ must be
associated with these gamma functions.
Let  $z \in \C^*$ and $[z] = z/\vert z \vert,~\vert z \vert_{\C} = \vert z \vert^2$.  Then the characters of $ \C^* $ have the form
$$
\chi : z \mapsto [z]^l\vert z \vert^t_{\C}   \quad  \quad   l \in \Z   \quad  t \in \C .
$$
As a preliminary definition  of the local $ L $-function corresponding to
unramified extensions at the point $ P $ of  $ X $, one can offer the same
gamma-factors for the character $ \chi $ as in the one-dimensional situation of the
local field $ \C $ \cite {Kn}:
$$
L_P(s, \chi) = 2(2\pi)^{-(s + t + \vert l \vert /2)}\Gamma (s + t + \vert l \vert /2).
$$

Here we assume that the unramified "extensions" of the field $ K_P = \C((t_P)) $
can be described by the Weil group of its residue field $ \C $. As we noted above,
although the field $ \C $ has no nontrivial extensions, the role of the Galois group played by its Weil group when we
build the local Langlands correspondence.

Next, we can assume that global $ L $-function
 of a sheaf  ${\cal F}$ on $ X $ will decompose into a product
over  all points of the curve of  these local $ L $-functions.
\section{Direct image conjecture}
%\begin{center}
%\vskip 0,5cm
%{\large {\bf Direct image conjecture}}
%\vskip 0,5cm
%\end{center}
Let $ f: X \rightarrow B $  be a proper morphism of  a surface  $ X $ onto a regular curve $ B $ with a smooth general fiber.
Return back to the hypothesis expressed above about the existence of direct
image from the set of abelian characters of $ K_2 (\A_X) $ to
the set of automorphic representations of $ \GL(\A_B) $ on $ B $. Let us formulate a number of properties which must be satisfied. The our starting point
are   well-known properties of the direct image (= automorphic induction) in one-dimensional case. Here they are (see, for example, \cite{AC}, \cite{HH}[Prop. 4.5]).

So, let $ f: B' \rightarrow B $ be a  finite unramified covering of $ B' $ and $ B $ of  degree $ m $. Then the direct image  $ f_* $ (= automorphic induction $ AI $) is a mapping $ AI: {\cal AF}^{nr}(n, B') \rightarrow {\cal AF}^{nr}(nm, B) $, and the inverse image  $ f ^ * $ (= base change $ BC $) is a mapping  $ BC: {\cal AF}^{nr}(n, B) \rightarrow {\cal AF}^{nr}(n, B') $. Here, $ {\cal AF}^{nr}(n, B) $ is the  space of  all unramified automorphic forms on $ \GL(n, \A_ {B}) $:
$$
{\cal AF}^{nr}(n, B) = \{\text{smooth functions on }~\GL(n, {\cal O})\backslash \GL(n, \A_ {B})/\GL(n, K)\},
$$
where  ${\cal O} = \prod_x{\cal O}_x,~x \in B$ and $K = \F_q(B)$.
We have the following properties:
\begin{description}
\item[i)] if  $\phi \in {\cal AF}^{nr}(n, B')$ and  $\psi \in {\cal AF}^{nr}(k, B)$, then the projection formula holds$$
f_*(\phi \boxtimes f^*(\psi)) = f_*(\phi) \boxtimes \psi\footnote{The operation $\boxtimes$ creates from the forms on the groups $\GL(n,~\A_ {B})$ and
$\GL(n',~\A_ {B})$
a form on the group  $\GL(n, \A_ {B})\times\GL(n', \A_ {B})$ which corresponds to the external tensor product of representations.};
$$
\item[ii)]  $f_*({\check \phi}) = {\check f_*(\phi)}  \quad \quad f_*({\check \phi}) = {\check f_*(\phi)}$,  where  ${\check \phi}(g) =\phi (^{t}g^{-1})$ (duality formula);
\item[iii)] if $n = 1$, then  $f_*({\chi})(gz) = \omega(z) f_*(\chi)(g)$ and $f_*(1)(gz) = \omega_0(z) f_*(1)(g)$ where $z \in \A_B^* =$ center of the group $\GL(\A_B)$ and $\omega\omega_0^{-1} = \beta^*(\chi)$ for the natural embedding $\beta : \A_B^* \rightarrow \A_{B'}^*$;
\item[iv)] for  unramified coverings  $f': B'' \rightarrow B'$ and $f: B' \rightarrow B$, we obtain
$$
(f'\circ f)^* = f'^*\circ f^*  \quad  \quad   (f'\circ f)_* = f_*\circ f'_*;
$$
\item[v)] if $k: B' \rightarrow B'$ and  $j: B \rightarrow B$ are  automorphisms of curves for which  $f \circ k = j \circ f$, then the following  base change rules hold\footnote {Let us make a terminological remark. Here, the expression "base change" means as it is customary in algebraic geometry. The same expression  was  understood  above as it is usually done in  the theory of automorphic representations.}
$$
k^* \circ f_* = j^*\circ f_* \quad \quad \text{and} \quad \quad  k_* \circ f^* = j_*\circ f^*;
$$
\item[vi)]  if $\phi \in {\cal AF}^{nr}(n, B')$ then $L_{\GL(n,~\A_{B'})}(\phi) =  L_{\GL(nm,~\A_B)}(f_*(\phi))$.
\end{description}
For the nonramified case, the definition of $L$-functions is given below on the section on the Hasse-Weil conjecture. These properties are supposed to be true for much more general case of ramified coverings. As far as I know, they were established, in the number fields case, only for the Galois cyclic extensions \cite{AC}.  In the geometric case, the existence of automorphic induction with these properties   follows from the Lafforgue's theorem (see footnote 21) and the corresponding construction for the Galois groups.

We show how to obtain the property iii) from the corresponding property for the direct images of $l$-adic sheaves. We have the following equality well-known  in the theory of representations:
$$
\Det\rho(g)\Det r(g)^{-1} = \chi(\text{Ver}(g)) \eqno (8)
$$
 for representations $ \rho $ and $ r$ of a group $ G $. Here, $ H \subset G $ is a  subgroup of finite index, $ \chi: H \rightarrow \C^* $ is a  character, $ \rho = \text{ind}^G_H (\chi),
r = \text{ind }^G_H (1), g \in G $ and $ \text{Ver}: G^{ab} \rightarrow H^{ab} $ is the transfer (Verlagerung) (definition of the latter, see
\cite{AT}[ch. XIII.2]).

For a finite unramified covering $ f: B'\rightarrow B $, we have an embedding of the fundamental group $ \pi_1(B) \rightarrow \pi_1(B') $ to a subgroup of finite index.  Applying  (8) to   their representations  we get the following  formula for the determinant of the image of a locally constant $ l$-adic sheaf $ {\cal F} $ of rank 1 on the curve $ B'$
$$
\Det f_* {\cal F} \otimes (\Det f_* \Q_l)^{-1} = \Nm_{B'/B}{\cal F}, \eqno (9)
$$
where $ \Det ({\cal G})_b = \Det ({\cal G}_b), ~ b \in B $ and the norm of the sheaf is defined as
$$
(\Nm_{B'/B} {\cal F})_b = \otimes_{b' \mapsto b} \Nm_{k(b')/k(b)} {\cal F}_{b'}.
$$
To get from this the  property iii) one needs to apply  the commutative diagram arising from the class field theory,
$$
\xymatrix @ C+1pt{
\text{sheaves } {\cal F} \text{ of rank 1 on } B'\ar[r] \ar[d] & \text {characters } \chi \text{ of } K_1(\A_{B'} )
\ar[d] \\ \text{sheaves } \Nm_{B'/B} {\cal F} \text{ on }  B \ar[r] &
\text{characters } \chi \vert_{K_1(\A_B)} \text{ of } K_1(\A_B).
}, \eqno (10)
$$
where the horizontal arrows are abelian Langlands correspondences, and to build the right vertical arrow we made  use  the natural embedding of $ K_1 (\A_B) \rightarrow K_1(\A_{B'}) $.
 Then one  applies  the formula for the central character of automorphic representations of the adele group that corresponds  to the $ l$-adic sheaf according to the Langlands correspondence (Lafforgue’s theorem). That is, the central character is the image of the sheaf
 $ \Det f_* {\cal F} $ of rank 1 with respect to the abelian Langlands correspondence.
 \vskip 2mm
When we turn to two-dimensional case, the number of functorial properties  will be increasing and the  properties itself are more complicated. Recall that we consider  maps  $f^i_*,~ i = 0, 1, 2 $ from the group ${\cal AF}^{nr}(K_2(\A_X))$ of unramified characters of  $ K_2 (\A_X) $, satisfying the reciprocity law (2) (and trivial on subgroup $K_2(\A_{12, X})$), to  the space of unramified automorphic forms ${\cal AF}^{nr}(\GL(r_i, \A_B)$.  Note that  the ranks  $ r_i $ depend, in general, on $ \chi $. We denote by $g$ the genus of a general fiber of the map $f$.
\vskip 0,5cm
{\bf Direct image conjecture}. For smooth proper  morphisms $f: X \rightarrow B$ of a smooth surface $X$ onto a smooth proper curve $B$, defined over a finite field   $\F_q$, there exist homomorphisms  $f^i_* : {\cal AF}^{nr}(K_2(\A_X)) \rightarrow {\cal AF}^{nr}(\GL(r_i, \A_B)$  ($~ i = 0, 1, 2 $) such that
\begin{description}
\item[o)] for fiber-trivial $\chi$, $r_0 = r_2 = 1$ and $r_1 = 2g$. For non-trivial $\chi$, $r_0 = r_2 = 0$ and $r_1 = 2g - 2$;
\item[i)] if  $\chi \in {\cal AF}^{nr}(K_2(\A_X))$ and  $\eta \in {\cal AF}^{nr}(K_1(\A_B))$, then the projection formula has the form
$$
f^i_*(\chi \otimes f^*(\eta)) = f^i_*(\chi) \otimes (\eta \circ \text{det}),
$$
where  $\text{det} : \GL(r_i, \A_B) \rightarrow \A^*_B$ is the determinant;
\item[ii)] for the direct images of a  character  $\chi \in {\cal AF}^{nr}(K_2(\A_X))$, the duality formula holds
$$
{\check f}^i_*(\chi) = f^{2-i}_*(\chi^{-1}\otimes \vert \cdot \vert),
$$
where  $\vert \cdot \vert$ is the  canonical character  (norm) on the group  $K_2(\A_X)$ (see  \cite{P1});
\item[iii)] $f^i_*(\chi)(gz) = \omega_i(z)f^i_*(\chi)(g)$ and $f^i_*(1)(gz) = \omega^0_i(z)f^i_*(1)(g)$ where  $z \in \A_B^* =$ center of the group $\GL(\A_B)$ and
$$
\prod_i \omega_i^{(-1)^i}(z) (\prod_i \omega_i^{0~(-1)^i}(z))^{-1} = \chi((\omega_{X/B}), \beta^*(z)),
$$
  where $\beta : \A_B^* \rightarrow (\A_X^{02})^*$ is the natural embedding, $(\omega_{X/B}) \in (\A_X^{01})^*$ is an idele corresponding to the relative canonical class (see \cite{P1}) and  $(-,-): \A_X^* \times \A_X^* \rightarrow K_2(\A_X)$ denotes multiplication in the   $K$-theory;

\item[iv)] for a finite unramified covering   $g: B \rightarrow B'$ of degree  $m$ we have
$$
(f\circ g)^i_*(\chi) =  g_*\circ f^i_*(\chi)
$$
on the group  $\GL(r_im, \A_{B'})$\footnote{A composition rule for finite unramified coverings of the scheme $ X $ must also be the case, but it involves going beyond the abelian group
$K_2(\A_X)$.};
\item[v)] if $k: X \rightarrow X$ and   $j: B \rightarrow B$ are automorphisms of schemes for which  $f \circ k = j \circ f$, then the following formula for base change is true
$$
k^* \circ f^i_* = j^*\circ f^i_* \quad \quad \text{and} \quad \quad k_* \circ f^* = j_*\circ f^*;
$$
\item[vi)] \hspace{4 cm} $L_X(\chi) = \prod_i L_{\GL(r_i,~\A_B)}(f^i_*(\chi))^{(-1)^i}$;
\item[vii)] there must also be a simple formula for the behavior under   blowings  of the  scheme $ X $ (the reader, tinker !).
\end{description}
\vskip 4mm
In fact, this hypothesis can be proved. Namely, we have
\vskip 2mm
{\bf Theorem}. {\em The existence of maps $ f^i_* $, satisfying the properties 0) - vi) for smooth proper morphisms $ f: X \rightarrow B $,  follows from the existence of the global Langlands correspondence for algebraic curve $ B$ over finite fields}.
\vskip 2mm
We prefer to consider this statement as a hypothesis, as there is a very non-trivial task to get an explicit construction of the  direct image, without using the Langlands correspondence and the general theory of etale cohomology. We now give an outline of how to get this construction, using Langlands correspondence and its properties. Referring to the diagram (5), let us define the mappings  $ f^i_* (\chi) $, using the two horizontal maps in this diagram  and the left vertical map. Upper arrow, the two-dimensional abelian Langlands correspondence is a bijection, which allows to invert it and thus to begin the construction of  the direct image. Left arrow follows from the theory of \'etale $ l$-adic sheaves. Finally, the most non-trivial part of the construction  is the bottom arrow,  the global Langlands correspondence for $ B$ (Lafforgue’s theorem \cite {L} \footnote{In the Lafforgue’s work  the correspondence  was constructed for  irreducible representations of the Galois group and cuspidal  representations  of adele group  on the curve. Apparently, the parabolic induction technique allows us to reduce the general case of completely irreducible representations to this one. This follows from Yu. G. Zarhin's semi-simplicity theorem for the Tate module of abelian varieties.}).

We sketch proofs of the most non-trivial properties of the direct image.

 For the property o),  we say that the character $\chi$ is fiber-trivial if its restriction on a fiber $F_b$ is trivial on the group $\Pic^0(F_b)$ for some $b \in B$. The property follows easily from the known facts on the cohomology of  l-adic sheaves.

  Since $ X $ is a surface over a finite field,  the duality formula ii) follows from  the relative duality theorem for $ l $-adic sheaf ${\cal F}$ on the scheme $ X $, corresponding to the character $ \chi $:
$$
\begin{array}{l}
R^if_*({\cal F}) \otimes R^{2-i}f_*(\Hom({\cal F}, \Q_l)(1)) \rightarrow R^2f_*({\cal F}\otimes \Hom({\cal F}, \Q_l)(1)) \\
\rightarrow R^2f_*( \Q_l)(1)) = \Q_l
\end{array}
$$
and subsequent application of the Langlands correspondence.

    In order to obtain the property iii), it is necessary to formulate an analogue of (9) in the new situation. We set
$$
(\Nm_{X/B} {\cal F})_b = \otimes_{x \in f^{-1}(b)} {\cal F}_x^{\nu_x(\omega_{X/B})} ,
$$
where ${\cal F}$ is a sheaf of rank 1 on $X$, $ b \in B$, $ \omega_{X/B} $ is a  section of the relative cotangent bundle and $ \nu_x $ is the valuation  at the point $ x$ on a  fiber of the mapping $ f$. We then have
$$
\otimes_i(\Det R^if_*{\cal F})^{(-1)^i} \otimes_i(\Det R^if_*\Q_l)^{(-1)^{i+1}} = \Nm_{X/B}{\cal F}.
$$
This relation follows when we  apply to the fibers of the morphism $ f$  Deligne's formula for a locally constant sheaf $ {\cal F}$ of rank 1 on
a smooth projective geometrically irreducible curve $ C $ of genus $ g $ and defined over a field $ k $. It reads as
$$
\otimes_i(\Det \text{H}^i(C, {\cal F})^{(-1)^i} = \otimes_{x \in C}{\cal F}_x^{\nu_x(\omega)}(1-g), \eqno (11)
$$
where $ \omega $ is a non-zero differential form on $ C $ of degree 1 and $ {\cal F}(n) $ is  the Tate's twist of the  sheaf $ {\cal F} $ (see, e.g,, \cite{P1}[n 3.1]\footnote {
If $ k = \F_q $ and the sheaf  $ {\cal F} $ corresponds to a character $ \chi $ of the idele group $ K_1(\A_C) $, then (11) is a reformulation  of the expression for the elementary factor of the classical functional equation $ L_C (s, \chi) = \chi ((\omega))q^{s(2-2g)}q^{g-1}L_C(1-s, \chi^{-1})$ for $ L$-functions $ L_C(s, \chi)$ of the curve $ C $. The elementary factor $ \chi((\omega))q^{s (2-2g)}q^{g-1}$ corresponds to an action of the Frobenius automorphism of
the ground field on the one-dimensional $ l$-adic space dual to the space from the LHS in (11). The factor $ \chi((\omega))q^{s(2-2g)}$ is associated with the product
   $ \otimes_{x \in C}{\cal F}_x^{\nu_x(\omega)}$, and the twisting $(1-g)$ corresponds to $q^{g-1}$.}).

Morphism $ f $ gives the inverse map
$$
i: K_1(\A_B) \rightarrow K_2(\A_X), \quad \quad i(a) = (a, (\omega_{X/B})_ {02}), ~ a \in K_1(\A_B )
$$
where $ (\omega_{X/B})_{02}$ $ K_1(\A_{X, 02}) $ is the idele associated with the  section $ \omega_{X/B} $ (see \cite{P1}[n 2.2]).
We can now formulate the analogue of diagram (10)
$$
\xymatrix @ C +1 pt {
\text{sheaves } {\cal F} \text{ of rank 1 on} X \ar[r] \ar[d] & \text{characters } \chi \text{ of } K_2(\A_X)
\ar[d] \\ \text{sheaves } \Nm_{X/B} {\cal F} \text{ on } B \ar[r] &
\text{characters } i^* \chi \text{ of } K_1(\A_B).
}
$$
Then the same arguments as in the case of a finite covering, give the property iii).
\vskip 2mm
 To understand the property  vi), the group $K_2(\A_X)$ can be presented as  an adelic product over all flags $x \in C$ and then computed in the following way
 $$
% \prod_{x \in C}^{'}K_2(K_{x, C})  \rightarrow  \prod_{x \in C}^{'}K_2(K_{x, C})/K_2({\cal O}_{x, C})
 K_2(\A_X) \rightarrow  K_2(\A_X)/K_2(\A_{12})
    = \oplus_{x \in C} \Z
=  \oplus_{x \in X} \oplus_{C \ni x} \Z \rightarrow  \oplus_{x \in X} \Z\cdot 1_x.
 $$
 Here, the last arrow is the sum of all $x, C$-components around the given point $x$. Finally, since the character    $\chi$ is unramified and automorphic the reciprocity law implies that $\chi$ is defined on the latter group \footnote{and actually on its quotient $K_2(\A_X)/K_2(\A_{12})K_2(\A_{01})K_2(\A_{02})$  which is isomorphic to the second Chow group. This group is an extension of $\Z$ by a finite group according to the S. Bloch and K. Kato  finiteness theorems.}.

 Then,
 $$
 L_X (\chi): = \prod_{x \in \vert X \vert} (1 - \chi(1_x) )^{-1},
 $$
 where $1_x$ is a generator of the local group at the point $x \in X$.  The $L$-functions $L_{\GL(r_i,~\A_B)}(f^i_*(\chi))$ are the standard
$L$-functions attached to  the representations generated by the forms $f^i_*(\chi))$ (see \cite {GJ} and the next section).

This property follows from the analogous property
$$
L_X({\cal F}) = \prod_i L_B(R^if_*({\cal F}))^{(-1)^i}
$$
for  the corresponding $l$-adic sheaves ${\cal F}$ and $R^if_*({\cal F})$ and the Langlands correspondences on $X$ and $B$ which imply that $L_X({\cal F}) =  L_X (\chi)$ and
$L_B(R^if_*({\cal F})) =  L_{\GL(r_i,~\A_B)}(f^i_*(\chi))$ for all $i$.
\vskip3mm
{\bf Remark 5}. It is not difficult to state  local and semi-local versions of this conjecture (as we deal above) and formulate their interrelations.
We can even give an explicit expression for these direct images in the case of  a  (smooth) mapping $f: X_b \rightarrow \Spec{\cal O}_b$ with general fiber of genus $g $, and the character $\chi$ is trivial (or is the  $s$-th degree of the norm $\vert \cdot \vert$). Then,  $f^0_*(\chi) = \vert \cdot \vert^s,~f^1_*(\chi) = \vert \cdot \vert^{s-1}$, and  $f^1_*(\chi)$ is   an Eisenstein series generating an unramified principal series representation  of group  $\GL(2g, K_b)$. The Satake parameters of  this representation are  eigenvalues  of the Frobenius automorphism of the field $k(b)$, acting on  1-dimensional cohomology of the curve
$F_b = X_b \times_{\hat{\cal O}_b} k(b)$. Still, we can define them without any application of the  $l$-adic cohomology making use of  the Tate-Iwasawa
method on the curve $F_b$. The appropriate version of this method was given in  \cite{P4}.
%It gives an explicit construction for the direct image.
\vskip3mm
{\bf Remark 6}. Certainly, the choice of smooth morphisms in our conjecture is too restrictive (since we have a complete base $B$). There would be more reasonable to consider the semi-stable families $X$ over a proper curve $B$.  In this case, we can expect that the target of the direct image (still for unramified $\chi$ on $X$) will be tamely ramified automorphic representations $(\pi, V_{\pi} = \otimes_{b \in B}V_{\pi, b})$ of the groups $\GL(\A_B)$. If $\text{dim}_{\C}V_{\pi, b}^I = 1$ for the Iwahori subgroups  $I \subset \GL({\cal O}_b)$ and all the points $b$ of bad reduction then we attach to   $\pi$ a tamely ramified automorphic function on  $\GL(\A_B)$. Since for general tamely ramified representations $V_{\pi}$ of a group $\GL(K_b)$, $\text{dim}_{\C}V_{\pi}^I  >1$, to construct such function we need additional reasons. This is possible for the semi-stable  elliptic families  \cite{L2, G} where the special representations are the local components for the points of bad reduction.
\vskip3mm
{\bf Remark 7}. In the geometric case, the direct image conjecture follows from the two-dimensional class field theory and the Langlands correpondence on the curve $B$ (Lafforgue's theorem).  In general, if it were possible to determine not only the abelian automorphic representations, but also automorphic representations of the non-commutative groups of adeles on schemes $ X $ of any dimension (by developing, for example,  Kapranov's idea), then one could try to define the operation of the direct image  $ f_*^i $, expanding  the map $ f $ in a composition of the projection of a projective bundle and of an embedding. This is suggested by the Grothendieck's proof of the  Riemann-Roch theorem.
\section{A link with the Hasse-Weil conjecture}
% \begin{center}
%\vskip 1cm
%{\large {\bf A relation with the Hasse-Weil conjecture}}
%\vskip 0,5cm
%\end{center}
We assume that the conjecture on the direct image can be stated and proved in much more general situation. Namely, one can remove the smoothness condition for schemes and morphisms (see Remark 6) and thus consider the ramified characters. Finally, there  is an extension of the conjecture to the number field case.
Then, an (explicit) construction of direct images of this type would provide an opportunity to prove
the  Hasse-Weil  conjecture~\cite{W2} on meromorphic continuation and functional
equation for $L$-functions of algebraic curves defined over algebraic number fields
 (and reprove the corresponding Grothendieck theorem on
$L$-functions of algebraic curves defined over fields of algebraic
functions with a finite field of constants, {\bf without} any  use of  $l$-adic cohomology).

Curves with this field of definition have   models, two-dimensional schemes $X$ with a structure morphism on a one-dimensional
scheme $B$, either a curve over a finite field or a spectrum  of a  ring of integers in  an algebraic number field. The existence of the direct image
from $X$ to $B$ allows  to build an  automorphic form on $B$, starting from a character of the group $K_2(\A_X)$. The Mellin transformation
of this form  will be  the $L$-function of interest to us. Automorphy of the constructed in this way form shall, in accordance with the general theory, imply meromorphic continuation and functional equation of the $L$-function  \cite{W3, JL, GJ}.

We give a precise formulation of the corresponding theorem for the case of algebraic  curves defined over a finite field.
Let $ {\cal AF}(n, B) $ be the space of automorphic forms on $ \GL(n, \A_{B}) $, i.e. of the functions $ F $ on $ \GL(n, \A_{B}) $ such that:
\begin{itemize}
\item [i)] there exists a character $ \omega $ of the center $ Z  = \A_B^* $ of  $ \GL(n) $ such that $ F(zg) = \omega(z) F(g), ~ z \in Z,
~ g \in \GL(n, \A_{B}) $;
\item [ii)]  if $ \gamma \in \GL(n, K) $, where $ K = \F_q(B) $, then $ F(g\gamma) = F(g) $;
\item [iii)] there exists an open compact subgroup  $ \K'\subset \GL(n, ~ {\cal O}) $, where
 $ {\cal O} = \prod_x{\cal O}_x, ~ x \in B $, such that
$ F(kg) = F(g), ~ k \in \K', ~ g \in \GL(n, \A_{B}) $.
\end{itemize}
The group $ \GL(n, \A_{B}) $ acts on the space of automorphic forms by left translations.
Smooth representation of $ \GL(n, \A_{B}) $ is called {\ em automorphic} if it is embedded as a sub-representation in the space of automorphic forms $ {\cal AF}(n, B) $.

Each irreducible automorphic representation $ \pi $ is a tensor product $ \otimes_{x \in B} \pi_x $ of smooth irreducible representations of  local groups  $ \GL(n, K_x) $. For almost all $ x $ the representation $ \pi_x $ is a  spherical one \cite {Fl}, \cite {B} [Chap. 3.4].  According to the general theory, such a representation belongs to the principal series, i.e. can be obtained by parabolic induction from an  unramified character $ \chi: T \rightarrow \C^* $ of a maximal torus $ T $ of the group $ \GL(n, K_x) $ (see \cite{B} [theorem 4.6.4] for the case $ n = 2 $, the proof of which can be extended to an arbitrary $ n $). This character is given by a set of $ n $ complex numbers $ z_1, \dots, z_n $ (defined up to a permutation). Then we can define the local $ L $-function for spherical representations of $ \pi_x $ as
$$
L_B(s, \pi_x)  =  \prod_ i (1 - z_iq^{-s})^{-1}, ~ i  = 1, \dots,  n.
$$
and the global $ L$-function as the Euler  product
$$
L_B(s, \pi)  =  \prod_ {x} L_(s, \pi_x).  \eqno (12)
$$
Missing factors for a finite set of points $ x \in B $ are defined in a  more complicated way (see
\cite{JL} [Ch. 1, Theorem 2.18], \cite{G} [Ch. 2, Theorem 8.2 and Table A] for $ n = 2 $ and
\cite{GJ} [Ch. 1, Theorem 3.3] for an arbitrary $ n $).
\vskip 2mm
{\bf Theorem.} {\em Let $ B $ be a smooth projective curve over a finite field $ \F_q $ and let $ \pi $ be an automorphic representation of $ \GL(n, \A_B) $. Then:}
\begin{itemize}
\item [i)] {\em the Euler product (12) for the  $ L$-function $ L_B(s, \pi) $ converges in the right half-plane $ \text{Re} (s)  >  1 $};
\item [ii)] {\em  $ L_B(s, \pi) $ extends meromorphically to the whole $ s$-plane as a rational function of $ q^{-s} $};
\item [iii)] {\em $ L_B(s, \pi) = \epsilon (s, \pi) L_B (1 - s, \check{\pi}) $, where $ \epsilon (s, \pi) $ is an elementary factor of the form $ aq^{bs}, ~ a, b \in \R $}.
\end{itemize}
\vskip 2mm
The theorem was proven for cuspidal representations in \cite{GJ} [Theorem 8.13]. The general case reduces to this one by using  parabolic induction  \cite{J} [Theorem 6.2]. For the group $ \GL(2) $, the  result was obtained by Jacquet and Langlands in 1970 \cite{JL} [Ch. II, Theorem 11.1]. Note that in these works an analogue of this theorem for number fields was also  proven.

If  {\em all} the local components of the representation $ \pi $ are spherical, then it is generated by an unramified  automorphic form (=tensor product of the local spherical vectors) and its Mellin transform is equal to the $ L$-function $ L_B (s, \pi) $. This fact allows us to apply this theorem to automorphic forms  $ f^i_* (\chi) $, arising as direct images of an unramified automorphic character $ \chi $ of $ K_2(\A_X) $ on a surface $ X $, which are mapped onto  the curve $ B $ .

Such a construction can be seen as a way of solving   the general problem, discussed in \cite{P3}[5.2]. It is a generalization   of the Tate-Iwasawa approach \cite{T1, Iw} to the zeta- and $ L$-functions of one-dimensional schemes to the schemes of higher dimension (in this case of dimension 2). It seems that the harmonic analysis on two-dimensional adelic spaces  \cite{P3}[2] and the theory of representations of discrete Heisenberg groups \cite{P3}[3, 4] are the necessary tools to  build  this  (hypothetical) direct image.
\vskip 2mm
{\bf Remark 8}. Let  $X$    be an $n$-dimensional regular scheme and  let $f: X \rightarrow B$ be a plat proper morphism onto a one-dimensional regular scheme $B$ with smooth general fiber. In the number field case, the schemes $X$ and $B$ must be considered in the Arakelov geometry, completed by archimedean fibers .  The $n$-dimensional class field theory suggests that there exists a reciprocity map
 $$
\varphi_X : K_n^M(\A_X) \longrightarrow G_K^{ab},
$$
where  $K_n^M(\A_X)$ is the Milnor part of the    $n$-th  Quillen's K-functor and  $G_K^{ab}$ is the Galois group of the maximal abelian covering of the scheme  $X$. We may assume that the direct image conjecture is valid  (with the appropriate changes) also for these morphisms and automorphic characters of the group   $K_n^M(\A_X)$. Automorphy of the characters means that they are trivial on  $n$ subgroups
$$
K_n^M(\A_{0, \hat{1}, 2, \dots, n}),\quad  K_n^M(\A_{0,1, \hat{2}, \dots, n}),~\dots, \quad  K_n^M(\A_{0, 1, \dots, \hat{n}}),
$$
where  $\hat{i}$  means that the   $i$-th index must be deleted. These relations are a generalization of the formulas (2). Just as above, this conjecture  has to imply the Hasse-Weil conjecture for the scheme  $X$.

 The first step to get the conjecture is  to construct direct images $f_* \colon K_n(\A^{\bullet}_X) \to K_{n-k}(\A_B^{\bullet})[k]$ for  proper morphisms $f: X \rightarrow Y$ of relative dimension $k$,  full adelic complexes
and $n-k \geq 0$. We do not consider the reasonable conditions for existence of these direct images. Let us only mention that for  closed regular embeddings the map  $f_*$  depends on a choice of the local equations for the $Y$ in $X$. Nevertheless, the resulting homomorphism $f_*$  is defined uniquely up to a homotopy.
\section{Appendix: zero-dimensional generalization \\ of the Langlands correspondence}
%\begin{center}
%\vskip 3cm
%{\large {\bf  Appendix: zero-dimensional generalization of the Langlands correspondence}}
%\vskip 0,5cm
%\end{center}
Let us consider in more detail an analogue of the Langlands correspondence in  zero-dimensional situation, i.e. for finite fields and their extensions.
As in the case of fields of higher dimension class field theory gives us a canonical reciprocity map
$$
K_0(K) \longrightarrow G_K^{ab} =  G_K,\eqno (13)
$$
which explicitly looks as follows:
$$
\Z \longrightarrow {\hat \Z}
$$
and maps  $1$  into Frobenius automorphism $\text{Fr}~ (\text{Fr}(x) = x^q)$ for a field  $K = \F_q$.
Here  $K_0( \F_q)$ is the Grothendieck  $K$-group of category of finite dimensional vector spaces, and  $1 \in K_0( \F_q)$ corresponds to the trivial one-dimensional space.

If $ K'/ K $ is a  finite extension of degree $ m $, then we have the commutative diagrams
$$
\xymatrix@C-17pt{
 K_0(K') \ar[d]_{\text{Nm}}\ar[r] & G_{K'}\ar[d]_i & K_0(K') \ar[r] & G_{K'}\\
K_0(K) \ar[r] & G_{K}& K_0(K) \ar[r] \ar[u]_j & G_{K}\ar[u]_{\text{Ver}}
%\ar[l] \ar[d]_{f^i_*} & \text{ Характеры } K_2(\A_X) \ar[l] \ar[d]_{f^i_*} \\
%\parbox{3.7cm}{\text{Автоморфные} \text{\qquad формы } \text{\qquad на }~$\GL(r_i,K_b)$} &
%\parbox{3.8cm}{\text{ Автоморфные} \text{\qquad формы } \text{\qquad на }~$\GL(r_i,K_b)$} \ar[l]&
%\parbox{3.8cm}{\text{ Автоморфные} \text{\qquad формы } \text{\qquad на }~$\GL(r_i,\A_B)$.} \ar[l]
}
$$
where the embeddings  $i$ and $j$ induced by the inclusion of $ K $ in $ K '$, $ \text{Nm} $ is the  norm map and $ \text{Ver} $ is the  transfer from  $ G_K ^ {ab} = G_K $ to $ G_ {k '} ^ {ab} = G_ {K'} $. In this case,   $\text{Nm}$  is multiplication by $ m $, and $ j $ is isomorphism (of the  groups equal to $ \Z $).

As above, the map (13) leads to a homomorphism  of the groups of characters  (one-dimensional representations) in the opposite direction.
We can define functors of direct and inverse images  for the characters and   these diagrams  lead to their properties similar to those of the previous section.

As we remember, the Kapranov's proposal for the 0-dimensional correspondance was to take the $L$-function $L(s, \rho)$ of the Galois representation $\rho$. If  $\lambda_1, \dots, \lambda_n $ = spectrum  of $\rho$  then  $L(s, \rho) = \det (I - \rho(Fr)q^{-s})  = \prod_i(1 - \lambda_iq^{-s})$. On the other hand, there is a structure of $d$-representations ($d = 1$ or $2$) for the correspondances in dimensions $d$. We can expect that something like 0-representations might appear at the RHS of the 0-dimensional correspondance. To understand the position we give the necessary definitions of the $d$-representation theory ($d = 0, 1, 2$ ).

This is a pure algebraic theory that has nothing to do with arithmetic. We fix a group $G$, a ground field $k$ and a cardinal  number $n$. The inductive structure of the theory looks then as follows:
$$
\begin{tabular}{|c|c|c|c|c|} \hlx{hv}
$d$ & representations & $d$-representations & $n$ = dimension of & sets of  1-dim.  \\ \hlx{vv}
       & spaces             &                                 & representation   & representations    \\ \hlx{vhv}
$2$ & $\text{V-mod}/k$ & functor $\pi(g): $ & $ \text{rk}~\text{V-mod}/k$ & $H^2(G, k^*)$   \\ \hlx{vv}
&&   $\text{V-mod} \rightarrow \text{V-mod}$    && \\ \hlx{vhv}
$1$ & $V/k$  & homomorphism & $\text{dim}_k~V$ & $H^1(G, k^*)  =$  \\ \hlx{vv}
&& $\pi(g): V \rightarrow V$   &&$\Hom(G, k^*)$ \\ \hlx{vhv}
$0$ & $A \subset Sym^n(k^*)$ & just $A$, no $\pi(g)$ & $\text{deg}~ A = n$ &  $H^0(G, k^*) = k^*$ \\ \hlx{vhv}
\end{tabular}
$$
$$
\begin{tabular}{|c|c|c|c|} \hlx{hv}
$d$ & distinguished  & $\Hom'$s & sets of  \\ \hlx{vv}
    & object      &  & characters   \\ \hlx{vhv}
$2$   & Vect/k  & $\text{Func}({\cal C}, {\cal D})$ & $\{g \mapsto V_g \}$  \\ \hlx{vv}
$1$  & k & $\Hom(V, U)$ & $\text{Map}(G, k)$ \\ \hlx{vv}
$0$  & 1 & $(\sum_in_ia_i^{-1})(\sum_jm_jb_j)$  & $k$\\ \hlx{vhv}
\end{tabular}
$$
Let here  $\text{V-mod}/k$ be a category-module (see above),   $ \text{rk}~\text{V-mod}/k$ be its rank (see \cite{GK, K}),  $V_g$ be a vector space over $k$ depending on $g \in G$, $V/k$ or $V$ be a vector space over $k$ and $\text{Sym}^n(k^*)$  be the set of cycles of degree $n$ with non-negative multiplicities.
$\text{Func}({\cal C}, {\cal D})$ denotes the set of functors from the category ${\cal C}$ to the category ${\cal D}$ which commute with the action of the category  $\text{Vect}/k$ (in the categorical sense).

We see that for $d = 0$ the group $G$ itself doesn't participate in all the constructions. For $d = 0$, the direct sum of the representations $A$ and $B$  is the sum of cycles  $A$ and $B$, the tensor product is their Pontryagin product in the group $k^*$ and dual to $A = \sum_i n_i a_i $ is ${\check A} =  \sum_i n_i a_i^{-1}$. The trace $\text{Tr}(A)$ is $\sum_i n_i a_i$ where summation is done in the group $k$. For $d = 2$ the dual "object" to a category ${\cal C}$ is the category of functors $\Func({\cal C}, \text{Vect}/k)$ \footnote{I'm grateful to participants of my Steklov seminar (especially to Sergey Gorchinskiy and Denis Osipov) for a discussion that clarified some points in the final version of the table.}.

 Now, we try to apply this small theory to the local Langlands correspondence. First, the field $k = \C$. The group at the LHS of the correspondance is the same Galois group for all $d$. But at the RHS, we see 2-representations of the group $\GL(2n, K)$ for  $d = 2$ and 1-representations of the group $\GL(n, K)$ for  $d = 1$.  We  may  guess that for $d = 0$ one needs to consider 0-representations of the group $\GL(0, K)$, where the latter group is simply trivial.  Then the Langlands correspondence could look as follows:
 $$
\begin{array}{l}
\text{Rep}^{(n)}(G_K) = \{n\text{-dimensional semi-simple representations of } G_K \text{ over } \C\}\stackrel{?}{\Leftrightarrow} \\
\{\mbox{irreducible 0-representations of the group  }~ \GL(0, K) =(1)\},
\end{array}
$$
 where $\rho \in \text{Rep}^{(n)}(G_k)$  goes to the determinant $\lambda_1 \cdots \lambda_n $ of the spectrum $\lambda_1, \dots, \lambda_n $   of $\rho$.
We see that this construction  does not coincide  completely with the Kapranov' proposal.  The $L$-function from his choice doesn't belong to the field $\C$ as we could expect according to the table above. Our choice is the determinant  which is an element of $\C$. Also, it satisfies the same multiplicative property for exact sequences of representations  as the $L$-function.
But, the correspondence will not a bijection in contrast with the behavior  of the $L$-functions which are actually the characteristic polynomials of the Frobenius
automorphism.
This problem we leave as a question.

Let us say that we also meet  another difficulty.
The Langlands correspondance is in a perfect accordance with reciprocity law for $d = 1$ or $2$. However,  we do not see here the group $K_0(K)$ at the RHS of the correspondance for $d = 0$. This unexpected puzzle can be solved along the following line.

In his paper \cite{K}, Kapranov justifies the appearance of the non-commutative group  $\GL(2, K)$ at the RHS of the abelian (!) 2-dimensional correspondance by  the following relation
$$
K_2(K) = H_2(\GL(2, K), \Z)/H_2(\GL(1, K), \Z)
$$
which is a particular case of the general Suslin's theorem on Milnor's $K^M_n(K)$ for a field $K$. In case of $n = 2$ this could follow from the Matsumoto's theorem that $K^M_2(K) = K_2(K) = H_2(\SL(K), \Z)$. Applying this, Kapranov has shown that  1-dimensional 1-representations of the $K_2(K)$ are in one to one correspondance with   1-dimensional 2-representations of the group  $\GL(2, K)$ (with a trivial restriction  to the subgroup $\GL(1, K)$). Since  $K_1(K) = K^* = H_1(\GL(1, K), \Z)$  we have an analogous relation  for 1-representations for $d = 1$.

For $d = 0$ we can add  the relation
$$
K_0(K) = H_0(\GL(0, K), \Z) =  \Z
$$
and the 0-dimensional Langlands correspondance looks as follows
$$
\begin{array}{l}
\text{Rep}^{(n)}(G_K) = \{n\text{-dimensional semi-simple 1-representations of } G_K \text{ over } \C\}\Leftrightarrow \\
 \{n\text{-dimensional semi-simple 1-representations of } K_0(K)  \text{ over } \C\}\Rightarrow \\
\{\mbox{1-dimensional  0-representations of the group  }~ \GL(0, K) =(1)\}.
\end{array}
$$
The same versions (but with bijections) are valid also for $d = 1$ or $2$. Besides that, it is unusual that we have to consider all semi-simple
representations at the RHS  of the correspondence but not the irreducible ones as we have for the case of dimension 1.
%{\Remark 7}.   The spectral correspondence behaves very well with respect  to  finite extensions
 %$k'/k$ of degree $m$. Namely, an $n$-dimensional  0-representation $\lambda_1, \dots, \lambda_n $  for the  $\GL(0, k)$ goes to  $\lambda_1^m, \dots, %\lambda_n^m $ for the  $\GL(0, k')$ by the base change.  An $n$-dimensional  0-representation $\mu_1, \dots, \mu_n $  for the  $\GL(0, k')$ goes to  %$\mu_1\cdot \zeta_1, \dots, \mu_n\cdot \zeta_m $ by the automorphic induction.
%Here, $\zeta_i, ~i = 1, \dots, m$ runs through all $m$-th roots of unity.

%\begin{center}
%{\Large {\bf Литература}
%\vskip1cm}
%\end{center}

%\begin{center}

%\end{center}
\noindent Steklov Mathematical Institute,\\
Russian Academy of Sciences, \\
Russia 119991, Moscow, Gubkina str. 8\\
e-mail  $parshin@mi.ras.ru$
\end{document}